\documentstyle[12pt]{article}

       \font\tenmsb=msbm10
       \font\sevenmsb=msbm7
       \font\fivemsb=msbm5
       \catcode`\@=11
       \ifx\amstexloaded@\relax\catcode`\@=\active
       \endinput\else\let\amstexloaded@\relax\fi
       \def\spaces@{\space\space\space\space\space}
       \def\spaces@@{\spaces@\spaces@\spaces@\spaces@\spaces@}
       \def\space@.  {\futurelet\space@\relax}
       \space@.   %
       \def\Err@#1{\errhelp\defaulthelp@\errmessage{AmS-TeX error: #1}}
       \def\relaxnext@{\let\next\relax}
       \def\accentfam@{7}
       \def\noaccents@{\def\accentfam@{0}}
       \def\Cal{\relaxnext@\ifmmode\let\next\Cal@\else
       \def\next{\Err@{Use \string\Cal\space only in math mode}}\fi\next}
       \def\Cal@#1{{\Cal@@{#1}}}
       \def\Cal@@#1{\noaccents@\fam\tw@#1}
       \def\Bbb{\relaxnext@\ifmmode\let\next\Bbb@\else
       \def\next{\Err@{Use \string\Bbb\space only in math mode}}\fi\next}
       \def\Bbb@#1{{\Bbb@@{#1}}}
       \def\Bbb@@#1{\noaccents@\fam\msbfam#1}
       \newfam\msbfam
       \textfont\msbfam=\tenmsb
       \scriptfont\msbfam=\sevenmsb
       \scriptscriptfont\msbfam=\fivemsb

\newtheorem{Theorem}{Theorem}[section]
\newtheorem{Lemma}{Lemma}[section]

\newtheorem{Corollary}{Corollary}[section]
\newtheorem{Remark}{Remark}[section]

\@addtoreset{equation}{section}

\begin{document}
\setlength{\columnsep}{5pt}
\title{On the generalized resolvent of linear pencils in Banach spaces.}
\author{\small  Qianglian Huang,$^*$\ \  Shuangyun Gao}
\date{ }
\maketitle\footnotetext{\scriptsize{\bf $^{*}$ Corresponding author.
}}
 \maketitle\footnotetext{\scriptsize{\bf
This research is supported by the Natural Science Foundation of
China (10971182), the Natural Science Foundation of Jiangsu Province
(BK2010309 and BK2009179), the Tianyuan Youth Foundation (11026115),
the Natural Science Foundation of Jiangsu Education Committee
(07KJB110131 and 10KJB110012) and the Natural Science Foundation of
Yangzhou University.}}

\par\vspace{0.1in}
 \noindent{\bf Abstract}
  {\small \  Utilizing the stability
characterizations of  generalized inverses of linear operator, we
investigate the existence of generalized resolvents of  linear
pencils in Banach spaces. Some practical criterions for the
 existence of
 generalized resolvents of the linear pencil $\lambda\rightarrow
T-\lambda
 S$ are provided and  an explicit expression of the generalized resolvent
 is given. As applications, the characterization for
the Moore-Penrose inverse of the linear pencil to be its generalized
 resolvent and the
existence of the
 generalized resolvents of linear pencils of finite rank operators, Fredholm operators and
semi-Fredholm operators are also considered. The results obtained in
this paper extend and improve many results in this area. }

\par\vspace{0.1in}\noindent{\bf Keywords}
 generalized inverse, generalized resolvent, linear
pencils, Moore-Penrose inverse, Fredholm operator, semi-Fredholm
operator
  \par\vspace{0.1in}
 \noindent{\bf{2000 Mathematics Subject Classification}}\ 47A10,
 47A55.
\par\vspace{0.1in}
\section{Introduction and Preliminaries}
\quad\  Let $X$ and $Y$ be two Banach spaces. Let $B(X, Y)$ denote
the Banach space of all bounded linear operators from $X$ into $Y$.
We write $B(X)$ as $B(X, X)$. The identity operator will be denoted
by $I$. For any $T\in B(X,Y)$, we denote by $N(T)$ and $R(T)$ the
null space and the range of $T$, respectively. \par The resolvent
set $\rho(T)$ of $T\in B(X)$ is, by definition, the set of all
complex number $\lambda\in C$ such that  $T-\lambda I$ is invertible
in $B(X)$. And its resolvent $R(\lambda)=(T-\lambda I)^{-1}$ is an
analytic function on $\rho(T)$ since it satisfies the resolvent
identity:
$$R(\lambda)-R(\mu)=(\lambda-\mu)R(\lambda)R(\mu), \quad  \forall\lambda, \mu\in\rho(T).$$
 The spectrum $\sigma(T)$ is the complement of $\rho(T)$ in $C$. As
we all know, the spectral theory plays a fundamental role in
functional analysis. If the operator $T-\lambda I$ has a generalized
inverse, we can consider the generalized resolvent and generalized
spectrum. Some properties  of the classical spectrum $\sigma(T)$
remain  true in  the case  of the generalized one[1-4]. Recall that
an operator $S\in B(Y,X)$ is said to be an inner inverse of $T\in
B(X,Y)$ if $TST=T$ and an outer inverse if $STS=S$. If $S$ is both
an inner inverse and outer inverse of $T$, then $S$ is called a
generalized inverse of $T$[5]. We always write the generalized
inverse of $T$ by $T^+$. If $T$ has a bounded generalized inverse
$T^+$, then  $TT^+$ and $T^+T$ are projectors with $R(T T^+
 )=R(T)$, $R(T^+T)=R(T^+)$, $N(T^+T)=N(T)$, $N(TT^+)=N(T^+)$ and
 $$X=N(T)\oplus R(T^+), \quad Y=N(T^+)\oplus R(T).$$
\par Let's recall the concept of generalized resolvent. Let $T\in B(X)$ and $U$ be an open set in the complex plane. The
function
$$U\ni\lambda\rightarrow R_g(T, \lambda)\in B(X)$$ is said to be a
generalized resolvent of $T-\lambda I$ on $U$ if
\par $(1)$
$(T-\lambda I)R_g(T, \lambda)(T-\lambda I)=T-\lambda I$  for all
$\lambda\in U$; \par $(2)$ $R_g(T, \lambda)(T-\lambda I)R_g(T,
\lambda)=R_g(T, \lambda)$ for all $\lambda\in U$;\par $(3)$ $R_g(T,
\lambda)-R_g(T, \mu)=(\lambda-\mu)R_g(T, \lambda)R_g(T, \mu)$ for
all $\lambda$ and $\mu$ in $U$. \par\vspace{0.1in} The first two
conditions say that $R_g(T, \lambda)$ is a generalized inverse of
$T-\lambda I$ for each $\lambda\in U$, while the third one is an
analogue of the classical resolvent identity. We also refer to it as
the generalized resolvent identity, which plays an important role in
the spectrum since it assures that $R_g(T, \lambda)$ is locally
analytic. The generalized resolvents has been widely used in many
fields such as spectrum theory and theory of Fredholm
operators[1-4,6]. According to M. A. Shubin[7], there exists a
continuous generalized inverse function (satisfying $(1)$ and $(2)$
but not possibly $(3)$) meromorphic in the Fredholm domain
$\rho_\phi(T)=\{\lambda\in C: T-\lambda I$ is Fredholm$\}$. And he
points out that it remains an open problem whether or not this can
be done while also satisfying $(3)$, i.e., it is not known whether
generalized resolvents  always exist. Many authors have been
interested in the existence problem for the generalized resolvents
of linear operators in [1-4,6]. In [2], C. Badea and M. Mbekhta
  proved that $T-\lambda I$ has
an analytic generalized resolvent in a neighborhood of $0$ if and
only if $T$ has a generalized inverse and $N(T)\subset R(T^m),
\forall m\in N$. It is worth mentioning that the condition
$N(T)\subset R(T^m)$ is not easy to be verified and its geometric
significance is vague.  In [3], C. Badea and M. Mbekhta introduced
the concept of linear pencil and its generalized resolvents.
\par Let $T, S\in B(X,Y)$  and $U$ be an open set in the complex plane. The
function
$$U\ni\lambda\rightarrow G(\lambda)\in B(Y,X)$$ is called a
generalized resolvent on $U$ of the linear pencil
$\lambda\rightarrow T-\lambda S$ if
\par $(1)$
$(T-\lambda S)G(\lambda)(T-\lambda S)=T-\lambda S$  for all
$\lambda\in U$; \par $(2)$ $G(\lambda)(T-\lambda
S)G(\lambda)=G(\lambda)$ for all $\lambda\in U$;\par $(3)$
$G(\lambda)-G(\mu)=(\lambda-\mu)G(\lambda)SG(\mu)$ for all $\lambda$
and $\mu$ in $U$.
\par If $X=Y$ and $S=I$, then the two concepts are coincident. We
would like to remark that this is not just the formal extension
since $N(T)\subset R(T^m)$ does not hold in general if $X\neq Y$.
The linear pencils are very useful in the study of stability radius
of Fredholm operators[3]. Utilizing the reduced minimum modulus and
the gap between closed subspaces, C. Badea and M. Mbekhta proved the
following theorem:
\begin{Theorem}\label{The1.1}$([3])$
 Let $X$ and $Y$ be two Banach spaces. Let $T,S\in B(X,Y)$ and $U\subset C$ be an open
 set. There exists a generalized resolvent for $\lambda\rightarrow T-\lambda
 S$ on $U$ if and only if the linear pencil $\lambda\rightarrow T-\lambda
 S$ has fixed complements on $U$, i.e., there exist two closed
 subspaces $E$ and $F$ of $X$ and $Y$ such that for all $\lambda\in
 U$,
 $$X=N(T-\lambda S)\oplus E,\quad\ Y=R(T-\lambda S)\oplus F.$$
\end{Theorem}
\par
In this paper, we utilize the stability characterizations of
generalized inverses of linear operator to investigate the existence
of generalized resolvents of  linear pencil $\lambda\rightarrow
T-\lambda
 S$ in Banach spaces. Some
practical criterions for the
 existence of
 generalized resolvents of the linear pencils are provided and  an explicit expression of the generalized resolvent
 is given. As applications, the characterization for
the Moore-Penrose inverse of the linear pencil  to be its
generalized
 resolvent and the
existence of the
 generalized resolvents of linear pencils of finite rank operators, Fredholm operators and
semi-Fredholm operators are also considered. The results obtained in
this paper extend and improve many results in this area.
\section{Main Results}
\quad\ We start our investigation with the following lemmas, which
are preparation for the proofs of our main results.
\begin{Lemma}\label{Lem2.1}$([8-10])$ \ Let
$T\in B(X,Y)$ with two bounded generalized inverses $T^+$ and
$T^\oplus$. Then there exists a $\delta>0$ such that the following
 statements are equivalent$:$
\par\noindent\quad $(1)\ R(\overline{T})\cap N(T^+)=\{0\};$
 \par\noindent\quad $(2)\ R(\overline{T})\cap N(T^\oplus)=\{0\}$, \\
 where $\overline{T}\in B(X,Y)$
 satisfies $\| \overline{T}-T\|<\delta$.  \end{Lemma}
\begin{Lemma}\label{Lem1.2} $([8,11,12])$ Let $T\in B(X,Y)$
with a bounded generalized inverse $T^+$ and $\overline{T}\in
B(X,Y)$ with $\|T^+\|\|\overline{T}-T\|<1$. Then the following
statements are equivalent:\par\noindent\quad $(1)\
B=T^+[I+(\overline{T}-T) T^+]^{-1}=[I+T^+(\overline{T}-T)]^{-1}T^+$
is a generalized inverse of $\overline{T};$\par\noindent\quad $(2)\
R(\overline{T})\cap N(T^+ )=\{0\};$\par\noindent\quad $(3)\
Y=R(\overline{T})\oplus N(T^+);$
\par\noindent\quad $(4)$\ $X=N(\overline{T})\oplus R(T^+).$
\end{Lemma}
 \quad The following theorem not only provides a practical criterion for the
 existence of
 generalized resolvents for the linear pencils, but also
 gives an explicit expression of the generalized resolvent.
  \begin{Theorem}\label{The3.1}
  Let $X$, $Y$ be two Banach spaces and $T,S\in B(X,Y)$. \par\noindent\quad $(1)$ If the linear pencil $\lambda\rightarrow T-\lambda
 S$  has an
analytic generalized
 resolvent on a neighborhood of 0, then
 for any  generalized inverse $T^+$ of $T$, there exists a
neighborhood $U(0)$
 of  $0$ such that
 $$\ R(T-\lambda S)\cap N(T^+)=\{0\}, \quad \forall\  \lambda\in U(0); $$
 \quad $(2)$ If\, $T$ has a  generalized inverse $T^+$ and  there exists a
neighborhood $U$
 of  $0$ such that
 $$\ R(T-\lambda S)\cap N(T^+)=\{0\}, \quad \forall\  \lambda\in U, $$
 then the linear pencil $\lambda\rightarrow T-\lambda
 S$ has an analytic generalized
 resolvent on a neighborhood of 0.
 In fact, $$G(\lambda)=T^+(I-\lambda ST^+)^{-1}:\ Y\rightarrow
 X$$
 is a generalized resolvent of $\lambda\rightarrow T-\lambda
 S$ on a neighborhood of 0.
    \end{Theorem}
\quad {\bf Proof }(1).  If the linear pencil $\lambda\rightarrow
T-\lambda S$  has an analytic generalized
 resolvent $G(\lambda)$ on a neighborhood $V(0)$ of zero, then $R(T-\lambda S)\cap N(G(\lambda))=\{0\}$ and $G(0)$ is a generalized
 inverse of $T$, we denote it by $T^\oplus$.  From the  generalized resolvent
 identity:
 $$G(\lambda)-G(\mu)=(\lambda-\mu)G(\lambda)SG(\mu),\quad \forall\lambda,\mu\in  V(0) $$  we get
 $$N(G(\lambda))=N(G(\mu))\quad {\rm and}\quad R(G(\lambda))=R(G(\mu)),\quad \forall \lambda,\mu\in V(0).$$
 Particularly, $N(G(\lambda))=N(G(0))=N(T^\oplus)$, $\forall \lambda\in V(0).$ Hence for all $ \lambda\in V(0),$ $$R(T-\lambda S)\cap N(T^\oplus)=R(T-\lambda S)\cap N(G(\lambda))=\{0\}.$$
By Lemma 2.1, we can get
 for any  generalized inverse $T^+$ of $T$, there exists a
neighborhood $U(0)$
 of  $0$ such that
 $$\ R(T-\lambda S)\cap N(T^+)=\{0\}, \quad \forall\  \lambda\in U(0).$$
  \par\noindent\vspace{0.1in}\quad (2). If  $T$ has a  generalized inverse $T^+$ and there exists a neighborhood $U$
 of \ $0$ such that for all $\lambda\in U,$ $\ R(T-\lambda S)\cap N(T^+)=\{0\}.$
 Without loss of any generality, we can assume $|\lambda|<\|ST^+\|^{-1}, \forall\ \lambda\in U $. Then  by Theorem 2.2, for
each $\lambda\in U,$
$$G(\lambda)=T^+(I-\lambda ST^+)^{-1}:\ Y\rightarrow X$$
 is a generalized inverse of $T-\lambda S$ with $R(G(\lambda))=R(T^+)$ and $N(G(\lambda))=N(T^+)$.
 In the following, we shall
 show  that $G(\lambda)$
 is a generalized resolvent of the linear pencil $\lambda\rightarrow T-\lambda
 S$ on $U$. To this aim, it suffices
 to prove the  generalized resolvent identity, i.e., for all
$\lambda,\mu\in U$, $$G(
\lambda)-G(\mu)=(\lambda-\mu)G(\lambda)SG(\mu). $$  Set
$$P(\lambda)=(T-\lambda S)G( \lambda)\quad {\rm and}\quad
Q(\lambda)=G(\lambda)(T-\lambda S),$$  then $P(\lambda)$ is the
projector from $Y$ onto $R((T-\lambda S))$ with
$N(P(\lambda))=N(G(\lambda))=N(T^+)$ and $R(P(\lambda))=R(T-\lambda
S)$, and $Q(\lambda)$ is the projector from $X$ onto $R(T^+)$  with
$N(Q(\lambda))=N(T-\lambda S)$ and $R(Q(\lambda))=R(G(
\lambda))=R(T^+).$ Now we claim
$$P(\lambda)P(\mu)=P(\lambda)\quad {\rm and}\quad Q(\lambda)Q(\mu)=Q(\mu),\quad \forall \lambda,\mu\in U.$$
In fact, for any $y\in Y$,
$[P(\lambda)P(\mu)-P(\lambda)]y=-P(\lambda)[(I-P(\mu))y]$. Noting
$$[I-P(\mu)]y\in N(P(\mu))=N(T^+)=N(P(\lambda)),$$ we get
$P(\lambda)[(I-P(\mu))y]=0$. Then $P(\lambda)P(\mu)=P(\lambda)$. For
any $x\in X$, $Q(\mu)x\in R(G(
\mu))=R(T^+)=R(G(\lambda))=R(Q(\lambda))$, we get
$(I-Q(\lambda))Q(\mu)x=0$. Hence $Q(\lambda)Q(\mu)=Q(\mu)$. Thus
 \begin{eqnarray*}
(\lambda-\mu)G(\lambda)SG(\mu)&=&G(\lambda)[(T-\mu S)-(T-\lambda S)]G(\mu)\\
&=&G(\lambda)(T-\mu S)G(\mu)-G(\lambda)(T-\lambda S)G(\mu)\\
&=&G(\lambda)P(\mu)-Q(\lambda)G(\mu)\\
&=&G(\lambda)P(\lambda)P(\mu)-Q(\lambda)Q(\mu)G(\mu)\\
&=&G(\lambda)P(\lambda)-Q(\mu)G(\mu)\\
&=&G(\lambda)-G(\mu).
 \end{eqnarray*}
 Therefore, $$G(\lambda)=T^+(I-\lambda ST^+)^{-1}:\ Y\rightarrow
 X$$
 is a generalized resolvent of $\lambda\rightarrow T-\lambda
 S$ on a neighborhood $U$ of zero.   This
completes the proof.
 \par\vspace{0.1in} From Theorem 2.1 and Lemma 2.2, we can get the
 following corollary which is a generalization of Theorem 1.1 (i.e., Theorem 3.5 in [3]).
 \begin{Corollary}\label{Cor2.1}
  Let $X$, $Y$ be two Banach spaces and $T,S\in B(X,Y)$. Then
  the following statements are equivalent:
 \par\noindent\quad $(1)$\ the linear pencil $\lambda\rightarrow T-\lambda
 S$  has an
analytic generalized
 resolvent on a neighborhood of zero;
  \par\noindent\quad $(2)$\ $T$ has a  generalized inverse $T^+$ and there exists a
neighborhood $U$
 of  $0$ such that
 $$X=\ N(T-\lambda S)\oplus R(T^+), \quad \forall\  \lambda\in U;$$
   \par\noindent\quad $(3)$\ $T$  has a  generalized inverse $T^+$ and there exists a
neighborhood $U$
 of  $0$ such that
 $$Y=\ R(T-\lambda S)\oplus N(T^+), \quad \forall\  \lambda\in U;$$
\par\noindent\quad $(4)$\ for any  generalized inverse $T^+$ of $T$, there exists a
neighborhood $U$
 of  $0$ such that
 $$X=\ N(T-\lambda S)\oplus R(T^+), \quad \forall\  \lambda\in U;$$
   \par\noindent\quad $(5)$\ for any  generalized inverse $T^+$ of $T$, there exists a
neighborhood $U$
 of  $0$ such that
 $$Y=\ R(T-\lambda S)\oplus N(T^+), \quad \forall\  \lambda\in U;$$
 \par\noindent\quad  In this case, $$G(\lambda)=T^+(I-\lambda ST^+)^{-1}:\ Y\rightarrow
 X$$
 is a generalized resolvent of $\lambda\rightarrow T-\lambda
 S$ on a neighborhood of 0.
    \end{Corollary}
\begin{Remark}\label{Rem2.1} In [3], C.
Badea and M. Mbekhta proved that both $N(T-\lambda S)$ and
$R(T-\lambda S)$ have the fixed complement if and only if the linear
pencil $\lambda\rightarrow T-\lambda
 S$ has an analytic
generalized
 resolvent on a neighborhood of zero. It should be noted
that Corollary 2.1 shows that each one of $N(T-\lambda S)$ and
$R(T-\lambda S)$ has the fixed complement if and only if
$\lambda\rightarrow T-\lambda
 S$ has an analytic
generalized
 resolvent.
\end{Remark}
\begin{Theorem}\label{The3.2}
  Let  $T,S\in B(X,Y)$. Then there is an
analytic generalized
 resolvent for the linear pencil $\lambda\rightarrow T-\lambda
 S$ on a neighborhood of zero if and only if there exists a
neighborhood $U(0)$
 of  $0$ such that for all $\lambda\in  U(0)$,
 $T-\lambda
 S$ has the generalized inverse $(T-\lambda S)^+$ satisfying $$\lim_{\lambda\rightarrow0}(T-\lambda
 S)^+=T^+.$$
    \end{Theorem}
    \quad {\bf Proof } It suffices to prove the sufficiency. If there exists a
neighborhood $U(0)$
 of  $0$ such that for all $\lambda\in  U(0)$,
 $T-\lambda
 S$ has the generalized inverse $(T-\lambda S)^+$ satisfying $\lim\limits_{\lambda\rightarrow0}(T-\lambda
 S)^+=T^+,$
  then we put $$P_\lambda=I-(T-\lambda
 S)^+(T-\lambda
 S).$$
 Hence $\lim\limits_{\lambda\rightarrow0}P_\lambda=P_0.$
 Without loss of generality, we can suppose that for all  $\lambda\in  U(0)$, $\|P_\lambda-P_0\|\cdot\|P_0\|<1$
and $\lambda\cdot\|T^+\|\cdot\|S\|<1$. Then
\begin{eqnarray*}
P_0R(P_\lambda)&=&(I-T^+T)N(T-\lambda S)\\
&=&(I-T^+T+\lambda T^+S-\lambda T^+S)N(T-\lambda S)\\
&=&[I-T^+(T-\lambda S)-\lambda T^+S]N(T-\lambda S)\\
&=&(I-\lambda T^+S)N(T-\lambda S)
 \end{eqnarray*} and by the Banach Lemma,  the  operator
$$W=I-P_0+P_\lambda P_0=I+(P_\lambda-P_0)P_0$$
is invertible and its inverse $W^{-1}: X\rightarrow X$ is bounded.
Next, by Theorem 2.1,  we only need to show $$R(T-\lambda S)\cap
N(T^+)=\{0\},\quad\ \forall \lambda\in U(0).$$ Take $y\in
R(T-\lambda S)\cap N(T^+)$, then $y=(T-\lambda S)x$ and $
T^+(T-\lambda S)x=0,$ where $ x\in X$. Hence
\begin{eqnarray*}
&&T(I-\lambda T^+S)x\\
&=&T[I+T^+(T-\lambda S)-T^+T]x\\
&=&Tx+T^+(T-\lambda S)x-TT^+Tx=0
 \end{eqnarray*}
 which implies $$(I-\lambda T^+S)x\in N(T)=R(P_0).$$
Therefore
\begin{eqnarray*}
&&(I-\lambda T^+S)x\\
&=&P_0(I-\lambda T^+S)x\\
&=&P_0WW^{-1}(I-\lambda T^+S)x
\\
&=&P_0(I-P_0+P_\lambda P_0)W^{-1}(I-\lambda T^+S)x\\
&=&P_0P_\lambda P_0W^{-1}(I-\lambda T^+S)x\\&\in&
P_0R(P_\lambda)\\&=&(I-\lambda T^+S)N(T-\lambda S).
 \end{eqnarray*}
 By the invertibility of $(I-\lambda T^+S)$, we get $x\in N(T-\lambda
 S)$. Thus $y=(T-\lambda S)x=0$. This completes the proof.
 \begin{Remark}\label{Rem2.2} According to M. A. Shubin[7], there exists a continuous generalized
inverse function but not an analytic generalized
 resolvent. From Theorem 2.2, we can see that if there exists a continuous generalized
inverse function,  then we can find a relevant analytic generalized
 resolvent.
\end{Remark}
\quad Next we shall give the characterizations for the existence of
generalized  resolvents of the finite rank operators, Fredholm
operators and semi-Fredholm operators. Their proofs come directly
from Theorem 3.1 and Theorems in [8,11,13,14], we omit them.
 \begin{Theorem}\label{The3.2}
  Let  $T,S\in B(X,Y)$ and $T$ be a finite rank operator. Then there is an
analytic generalized
 resolvent for the linear pencil $\lambda\rightarrow T-\lambda
 S$ on a neighborhood of zero if and only if there exists a
neighborhood $U$
 of  $0$ such that
 $$\ {\rm Rank}(T-\lambda S)={\rm Rank}\ T, \quad \forall\  \lambda\in U. $$
    \end{Theorem}
\begin{Theorem}\label{The3.3}
  Let  $T,S\in B(X,Y)$ and $T$ be an Fredholm operator. Then there is an
analytic generalized
 resolvent for the linear pencil $\lambda\rightarrow T-\lambda
 S$ on a neighborhood of zero if and only if there exists a
neighborhood $U$
 of  $0$ such that for all $\lambda\in U$, either
 $$\ \dim N(T-\lambda S)=\dim N(T)\quad   or \quad {\rm codim} R(T-\lambda S)={\rm codim} R(T). $$
    \end{Theorem}
     \begin{Remark}\label{Rem2.2} Theorem 2.4 is a generaliztion of Theorem 4.1 in [3].
\end{Remark}
    \begin{Theorem}\label{The3.4}
   Let  $T,S\in B(X,Y)$ and $T$ be an semi-Fredholm operator with a generalized inverse. Then there is an
analytic generalized
 resolvent for the linear pencil $\lambda\rightarrow T-\lambda
 S$ on a neighborhood of zero if and only if there exists a
neighborhood $U$
 of  $0$ such that for all $\lambda\in U$, either
 $$\ \dim N(T-\lambda S)=\dim N(T)<\infty \quad or\quad  {\rm codim} R(T-\lambda S)={\rm codim} R(T)<\infty. $$
    \end{Theorem}
    \quad In the following,  we shall give the characterization for the
    Moore-Penrose inverse of the linear pencil $\lambda\rightarrow T-\lambda
 S$ to be its generalized
 resolvent. We recall that if the operator $T^\dag\in B(Y,X)$ satisfies
$$TT^\dag T=T,\quad T^\dag TT^\dag =T^\dag,\quad (TT^\dag)^*=TT^\dag,\quad
(T^\dag T)^*=T^\dag T,$$ where $T^*$ denotes the adjoint operator of
$T$, then  $T^\dag$ is called to be the Moore-Penrose inverse of
$T$. If $T^\dag$ is the Moore-Penrose inverse of $T$, then $T^\dag
T=P^\perp_{R(T^\dag)}$ and  $TT^\dag=P^\perp_{R(T)},$ where
$P^\perp_M$ is the orthogonal projector on $M$.
 \begin{Theorem}\label{The3.4}
   Let $X$ and $Y$ be two Hilbert spaces. Let $T,S\in B(X,Y)$ and $R(T)$ be closed. Then the Moore-Penrose inverse
 $(T-\lambda S)^\dagger$ of the linear pencil $\lambda\rightarrow T-\lambda
 S$ is its analytic generalized
 resolvent on a neighborhood of 0 if and only if there exists a
neighborhood $U$
 of  $0$ such that for all $\lambda\in U$,
 $$\  N(T-\lambda S)=N(T) \quad and\quad   R(T-\lambda S)=R(T). $$
    \end{Theorem}
    \quad {\bf Proof } If
    the Moore-Penrose inverse
 $(T-\lambda S)^\dagger$ of the linear pencil $\lambda\rightarrow T-\lambda
 S$ is its analytic generalized
 resolvent on a neighborhood $U$ of zero, then from the
 generalized resolvent identity,  we get
 $$N((T-\lambda S)^\dagger)=N(T^\dagger)\quad {\rm and}\quad R((T-\lambda S)^\dagger)=R(T^\dagger),\quad \forall \lambda\in U.$$
\noindent Hence $$\  N(T-\lambda S)=[R((T-\lambda
S)^\dagger)]^\perp=[R(T^\dagger)]^\perp=N(T)$$ and  $$R(T-\lambda
S)=[N((T-\lambda S)^\dagger)]^\perp=[N(T^\dagger)]^\perp=R(T).
$$
\indent Conversely, if there exists a neighborhood $U$
 of  $0$ such that for all $\lambda\in U$,
 $N(T-\lambda S)=N(T)$ and $R(T-\lambda S)=R(T),$ then $R(T-\lambda S)\cap N(T^\dagger)=R(T)\cap
 N(T^\dagger)=\{0\}$.
  By Theorem 2.1, $$G(\lambda)=T^\dagger(I-\lambda S T^\dagger)^{-1}=(I-\lambda T^\dagger
  S)^{-1}T^\dagger$$
is an analytic generalized
 resolvent of linear pencil $\lambda\rightarrow T-\lambda
 S$. Next we shall show that $G(\lambda)$ is also its Moore-Penrose
 inverse. Indeed, by $N(T-\lambda S)=N(T)$ and  $R(T-\lambda S)=R(T)$,
 then the orthogonal projector $P^\perp_{N(T-\lambda
 S)}=P^\perp_{N(T)}=I-T^\dagger T$ and $P^\perp_{R(T-\lambda
 S)}=P^\perp_{R(T)}=TT^\dagger$. Therefore,
  \begin{eqnarray*}
G(\lambda)&=&T^\dagger(I-\lambda S T^\dagger)^{-1}=T^\dagger
TT^\dagger(I-\lambda S T^\dagger)^{-1}=T^\dagger T(I-\lambda
T^\dagger
  S)^{-1}T^\dagger\\
&=&T^\dagger T(I-\lambda T^\dagger
  S)^{-1}T^\dagger TT^\dagger=[I-P^\perp_{N(T-\lambda
 S)}](I-\lambda T^\dagger
  S)^{-1}T^\dagger P^\perp_{R(T-\lambda
 S)}\\
&=&[I-P^\perp_{N(T-\lambda
 S)}]G(\lambda)P^\perp_{R(T-\lambda
 S)}.
 \end{eqnarray*}
So we get
\begin{eqnarray*}
G(\lambda)(T-\lambda S) &=&[I-P^\perp_{N(T-\lambda
 S)}]G(\lambda)P^\perp_{R(T-\lambda
 S)}(T-\lambda S)\\
 &=&[I-P^\perp_{N(T-\lambda
 S)}]G(\lambda)(T-\lambda S)\\
 &=&I-P^\perp_{N(T-\lambda
 S)}=T^\dagger T
 \end{eqnarray*}
 and
 \begin{eqnarray*}
(T-\lambda S)G(\lambda)&=&(T-\lambda S)[I-P^\perp_{N(T-\lambda
 S)}]G(\lambda)P^\perp_{R(T-\lambda
 S)}\\
 &=&(T-\lambda S)G(\lambda)P^\perp_{R(T-\lambda
 S)}\\
 &=&P^\perp_{R(T-\lambda
 S)}=TT^\dagger.
 \end{eqnarray*}
Hence $[G(\lambda)(T-\lambda S)]^*=G(\lambda)(T-\lambda S)$ and
$[(T-\lambda S)G(\lambda)]^*=(T-\lambda S)G(\lambda)$. Thus
 $G(\lambda)$ is  the Moore-Penrose
 inverse of  $ T-\lambda
 S$. This completes the proof.
 \begin{Corollary}\label{Cor2.2}
  Let $X$ be a Hilbert spaces and $T\in B(X)$, $R(T)$ be closed. Then the Moore-Penrose inverse
 $(T-\lambda I)^\dagger$ of the linear pencil $\lambda\rightarrow T-\lambda
 I$ is the analytic generalized
 resolvent on a neighborhood of zero if and only if
 $$\  N(T)=\{0\} \quad and\quad  R(T)=X. $$
 In this case, $T$ is invertible, the Moore-Penrose generalized
 inverse is the inverse and the  generalized
 resolvent is exactly its classical  resolvent.
    \end{Corollary}
      \quad {\bf Proof } If $N(T)=\{0\}$ and
      $R(T)=X$, by the inverse operator theorem, we can get what we
      desired. Conversely, if the Moore-Penrose inverse
 $(T-\lambda I)^\dagger$ is the generalized
 resolvent  of $T-\lambda
 I$ on a neighborhood of zero, then from Theorem 2.6, $N(T-\lambda I)=N(T)$ and $ R(T-\lambda
      I)=R(T)$ hold on a neighborhood of zero. It is easy to be verified that $N(T)=\{0\}$ and $R(T)=X$.  This completes the proof.
    \begin{Remark}\label{Rem2.2} Corollary 2.2 points out that why
    we use the generalized inverse instead of the Moore-Penrose inverse to
   define the generalized  resolvent. Just the nonuniqueness can give the generalized inverse more significance.
\end{Remark}
\subsection*{Acknowledgment}
\quad\ The authors would like to express their deep gratitude to
Professor Ma Jipu for his kind guidance.

 \par\vspace{0.2in}
 {\small
\noindent College of Mathematics\\
Yangzhou University\\
 Yangzhou 225002\\
 China
\\
E-mail: qlhmath@ yahoo.com.cn (Q. Huang),\\
\indent  \quad\quad\ gaoshuangyun@126.com (S. Gao).}

\begin{thebibliography}{BB}
 \bibitem{1} Mbekhta, M.,  {R${\rm\acute{e}}$solvant g${\rm \acute{e}}$n${\rm \acute{e}}$ralis${\rm \acute{e}}$ et th${\rm \acute{e}}$orie
spectrale,} J. Operator Theory. {21} (1989) 69-105.
\bibitem{2} Badea, C., Mbekhta, M., Generalized inverses and the maximal radius of regularity of a Fredholm
operator, Integr. Equ. Oper. Theory. {28} (1997), 133-146.
\bibitem{3} Badea, C., Mbekhta, M.,  {The stability
radius of Fredholm linear pencils,} J. Math. Anal. Appl. {260}
(2001) 159-172.
\bibitem{4} Mbekhta, M., {On the generalized resolvent in Banach spaces,} J. Math. Anal. Appl. {189} (1995), 362-377.
\bibitem{5} Nashed, M. Z.,  {Generalized Inverses and Applications,} Academic Press, New York, 1976.
 \bibitem{6} Hoefer, A., {Reduction of generalized resolvents of linear operator function,} Integr. Equ. Oper. Theory {48} (2004) 479-496.
\bibitem{7} Shubin, M. A., {On holomorphic families of subspaces,}
Integr. Equ. Oper. Theory. {2} (1979), 407-420.
\bibitem{8} Ma, J. P.,  {A rank theorem  of operators between Banach spaces,} Front. Math. China.
 {1} (2006) 138-143.
    \bibitem{9}  Huang, Q. L.,  Ma, J. P.,  {Continuity of generalized inverses of linear operators
in Banach spaces and its applications,}
   Appl. Math. Mech.  {26}(2005) 1657-1663.

\bibitem{10} Ma, J. P.,  { A generalized transversality in global analysis,}  Analysis in Theory and Applications.
{20}:4 (2004) 391-394.

\bibitem{11}Ma, J. P.,  {Complete rank theorem of advanced calculus and singularities of
bounded linear operators,} Front. Math. China. {2} (2008) 305-316.
\bibitem{12}Ma, J. P.,  {A generalized transversality in global analysis,} Pacific J. Math.
236:2 (2008) 357-371.

  \bibitem{13} Ma, Z. F., Ma, J. P., A common property of $R(E,F)$ and $B(R^n,R^m)$
 and a new method for seeking a path to connect two operators. Sci. China Math. 53:10 (2010), 2605-2620,

\bibitem{14}  Ma, Z. F., Ma, J. P., { The smooth Banach
submanifolds $B^*(E,F)$ in $B(E,F)$,}  Science in China. {52 A}
(2009) 2479-2492.







\end{thebibliography}
\end{document}